\documentclass[12pt]{article}
\usepackage[T1]{fontenc}
\usepackage[latin1]{inputenc}
\usepackage{lmodern,amsmath,amssymb}
\usepackage{graphicx}
\begin{document}

\title{On the plane-wave Riemann Problem \\ in Fluid Dynamics \footnote{Submitted for publication in 2009}}
\author{Bernd Einfeldt \\ b-einfeldt@t-online.de}
\maketitle  {\small \textbf{Abstract}: This paper contains a
stability analysis of the plane-wave Riemann problem for the
two-dimensional hyperbolic conservation laws for an ideal
compressible gas. It is proved that the contact discontinuity in the
plane-wave Riemann problem is unstable under perturbations. The
implications for Godunov´s method are discussed and it is shown that
numerical post shock noise can set of a contact instability. A
relation to carbuncle instabilities is established.

 $\\$ $\\$ \textbf{Key words:}
Riemann solver, Godunov-type methods, hyperbolic conservation laws,
gas dynamic, carbuncle instability $\\$ $\\$ \textbf{AMS(MOS)}
subject classifications: 65M08, 65M12
\newpage
\section{Introduction} Consider the compressible Navier-Stokes equations for a
viscous, heat conducting gas. The governing equations can be found
in many books. They are in two spatial dimensions
\begin{subequations}
\begin{align}
 &\text{Conservation of Mass}\nonumber \\
 &\frac{\partial}{\partial t}\rho + \frac{\partial}{\partial
 x}(\rho u) + \frac{\partial}{\partial y}(\rho v) = 0
 \label{eq:NV-1a} \\
 &\text{Conservation of x-Momentum}\nonumber \\
 & \frac{\partial}{\partial t}(\rho u) + \frac{\partial}{\partial
 x}(\rho u^2 + p)+  \frac{\partial}{\partial y}(\rho u v)
 = \label{eq:NV-1b} \\
 &\frac{\partial}{\partial x}[\frac{2}{3}\mu (2 \frac{\partial}{\partial x}u - \frac{\partial}{\partial y}v)]
 + \frac{\partial}{\partial y}[\mu (\frac{\partial}{\partial y}u + \frac{\partial}{\partial x}v)]
 \nonumber \\
 &\text{Conservation of y-Momentum} \nonumber \\
 &\frac{\partial}{\partial t}(\rho v) + \frac{\partial}{\partial x}(\rho u v) +
 \frac{\partial}{\partial y}(\rho v^2 + p)
 =  \label{eq:NV-1c} \\
 &\frac{\partial}{\partial x}[\mu (\frac{\partial}{\partial y}u + \frac{\partial}{\partial x}v)]
 + \frac{\partial}{\partial y}[\frac{2}{3}\mu (2 \frac{\partial}{\partial y}v - \frac{\partial}{\partial x}u)]
 &\nonumber \\
 &\text{Conservation of Energy} \nonumber \\
 & \frac{\partial}{\partial t}E + \frac{\partial}{\partial x}[u(E + p)] + \frac{\partial}{\partial y}[v(E + p)]
 = \label{eq:NV-1e} \\
 & \frac{\partial}{\partial x}[u \frac{2}{3}\mu (2 \frac{\partial}{\partial x}u - \frac{\partial}{\partial y}v)]
 + \frac{\partial}{\partial x}[v \mu (\frac{\partial}{\partial y}u + \frac{\partial}{\partial x}v)]
 \nonumber \\
 & \frac{\partial}{\partial y}[u \mu (\frac{\partial}{\partial y}u + \frac{\partial}{\partial x}v)]
 + \frac{\partial}{\partial y}[v \frac{2}{3}\mu (2\frac{\partial}{\partial y}v - \frac{\partial}{\partial x}u)]
 \nonumber \\
 &+ \frac{\partial}{\partial x}[k \frac{\partial}{\partial x}T]
  + \frac{\partial}{\partial y}[k \frac{\partial}{\partial y}T]
  \nonumber
\end{align} \label{eq:NV}
\end{subequations}
\newline
The dependent variable are the density $\rho$, the velocity field
$\vec{q}=(u,v)$ and the total Energy per unit volume E. The total
Energy per unit volume is given by
\begin{align}
    E=\rho e + \frac{1}{2}\rho(u^2 + v^2) \label{eq:E=rhop}
\end{align}
where $e$ is the internal energy per unit mass. The thermodynamic
variables $\rho$ and $e$ are related to the pressure $p$ through the
equation of state
\begin{align}
    p=p(\rho,e)\;\;\;\;\text{(equation of state)}
\end{align}
If we assume that the fluid is a perfect gas, the equation of state
state is
\begin{align}
    p=(\gamma-1)\rho e= (\gamma - 1)[E - \frac{1}{2}\rho (u^2 + v^2)] \;\;\;\;\mathrm{or}\;\;\;\; T=\frac{\gamma-1}{R}e  \label{eq:eqof}
\end{align}
where the constants $\gamma$ and $R$ are the ratio of specific heats
and the gas constant, respectively. T is the temperature. For a
perfect gas the specific heat at constant volume $c_v$ and the
specific heat at constant pressure $c_p$ are related to $\gamma$ and
$R$ by
\begin{align}
    c_v=\frac{R}{\gamma-1}\;\;\;\;\mathrm{and}\;\;\;\;c_p=\frac{\gamma R}{\gamma-1} \label{eq:gamme-R-relation}
\end{align}
We assume in the following that the coefficient of viscosity $\mu$
is constant. We also assumed that the coefficient of bulk viscosity
is negligible for the fluid, such that the second coefficient of
viscosity $\acute{\mu}$ is
\begin{align}
    \acute{\mu}=-\frac{2}{3}\mu
\end{align}
Furthermore we assume that Fourier's law for heat transfer holds,
with a constant coefficient of thermal conductivity $k$.
If the gas is polytropic resp. a perfect gas, the internal energy
$e$ is related to the temperature $T$ by $\eqref{eq:eqof}$; i.e.
\begin{align}
    e=c_v T    \label{eq:e-T-relation}
\end{align}
\\
For a vanishing viscosity $\mu$ and thermal conductivity $k$ the
Navier-Stokes equations reduce to the hyperbolic conservation laws
for an ideal compressible gas, also denoted as Euler equations. In
vector form these equations are
\begin{align}
   & \frac{\partial}{\partial t} \mathbf{u}(x,y,t) + \frac{\partial}{\partial x}\mathbf{f}(\mathbf{u}(x,y,t)) + \frac{\partial}{\partial y} \mathbf{g}(\mathbf{u}(x,y,t)) = 0\label{eq:Euler-diff}
\end{align}
where the conserved variable and flux functions are given by
\begin{align}
\mathbf{u}=\begin{pmatrix}
 \rho \\
 \rho u \\
 \rho v \\
 E
\end{pmatrix}
&\; \; \; \; \; \; \mathbf{f}(\mathbf{u})=\begin{pmatrix}
 \rho u\\
 \rho u^2 + p\\
 \rho uv\\
 u(E + p)
\end{pmatrix}
&\mathbf{g}(\mathbf{u})=\begin{pmatrix}
 \rho v\\
 \rho uv\\
 \rho v^2 + p\\
 v(E + p)
\end{pmatrix} \label{eq:con-var}
\end{align}
The velocity field will be denoted by $\mathbf{v}=(u,v)$. The
density $\rho$ and pressure $p$ are related to the conserved
quantities through the equation of state $\eqref{eq:eqof}$.
\newline
The Euler equations are a system of hyperbolic conservation laws;
i.e. for any value $\mathbf{u}_0=(\rho_0,\rho_0 u_0,\rho_0
v_0,E_0)^T$ with positive density $\rho_0$ and positive internal
energy $e_0$ the Jacobian matrix
\begin{align}
    \vec{n}D\mathbf{F}(\mathbf{u_0})=n_x D\mathbf{f}(\mathbf{u}_0) + n_y D\mathbf{g}(\mathbf{u}_0)
\end{align}
is diagonalizable with real eigenvalues for every unit vector
$\vec{n}=(n_x,n_y)$. Therefore the structure of a plane-wave
\begin{align}
    \mathbf{u}(x,y,t)=\mathbf{\varphi}(\vec{n}\cdot(x,y)-\dot{s}t)
\end{align}
propagating at speed $\dot{s}$, is independent of the orientation in
space; see for example [Lev2002]. The Euler equations are rotational
symmetric and Galilean invariant. A special class of plane-waves is
defined through the plane wave Riemann problem; an initial value
problem for $\eqref{eq:Euler-diff}$ with piecewise constant initial
data, separated by a straight line:
\begin{align}
    \mathbf{u}(\vec{x})=
    \begin{cases}
        \mathbf{u}_l & \text{for  } \vec{n}_0\cdot(\vec{x}-\vec{x}_0) < 0 \label{eq:per-rie-0}\\
        \mathbf{u}_r & \text{for  } 0 < \vec{n}_0\cdot(\vec{x}-\vec{x}_0)
    \end{cases}
\end{align}
where $\vec{x}_0$ is a given point in (x,y)-plane and $\vec{n}_0$ a
given unit-direction, $\mathbf{u}_l$ and $\mathbf{u}_r$ are initial
states at time $t=t^n$.
\newline
The integral form of the conservation law $\eqref{eq:Euler-diff}$ is
\begin{align}
    & \frac{\partial}{\partial t} \int_{D}^{}\mathbf{u}(\xi,\eta,t) \: d\eta \; d\xi  \label{eq:Euler-int} \\
    &= - \int_{\partial D}^{ } \vec{n}(\xi)\cdot \: \mathbf{F}(\mathbf{u}^R(\vec{x}(\xi),t+;\vec{n}(\xi))) \: d\xi \nonumber
\end{align}
where $\vec{n}(\xi)=(n_x(\xi),n_y(\xi))$ is the outward pointing
unit-normal vector on $\partial D$ at a point
$\vec{x}(\xi)=(x(\xi),y(\xi))$ on $\partial D$, where $\xi$ is the
arclength parametrization of the boundary $\partial D$. D is an
arbitrary bounded convex set with a piecewise smooth boundary.
$\mathbf{u}^R(\vec{x}(\xi),t+;\vec{n}(\xi))$ denotes the one sided
limit in time of the solution to the plane-wave Riemann problem
$\eqref{eq:per-rie-0}$ at the cell boundary $\vec{x}(\xi)$ in the
direction $\vec{n}(\xi)$; i.e.
\begin{align}
    \mathbf{u}^R(\vec{x}(\xi),t+;\vec{n}(\xi)):=\lim_{\epsilon\rightarrow0}\mathbf{u}^R(\vec{x}(\xi),t+\epsilon;\vec{n}(\xi)) \;
    \; \mathrm{\: on} \: \partial D
\end{align}
If the solution $\mathbf{u}$ is smooth, then we can apply the
divergence theorem and obtain the differential form
$\eqref{eq:Euler-diff}$ from the integral form
$\eqref{eq:Euler-int}$. If a solution of the Euler equations
contains discontinuities only the integral form is valid.
\newline
In this paper we consider the stability of the plane-wave Riemann
problem $\eqref{eq:per-rie-0}$ for the Euler equations and proof
that the plane-wave Riemann problem is unstable under perturbations.
A relation of this instability to the failing of Godunov´s and Roe´s
method reported by Quirk [QUI1994], XU [XU1999], Elling [VE2006],
Roe [ROE2007] is discussed. $\\$ The results of this paper can be
extended to general systems of hyperbolic conservation laws, which
have at least one genuinely nonlinear characteristic field and two
linear degenerate fields with a double eigenvalue. For clarity we
restrict the description to the Euler equation for a $\gamma$-law
gas $\eqref{eq:Euler-diff}$ which serves as a model equation for
more general systems. $\\$ \newline The outline of this paper is as
follows. In section 2 we analyze the stability of the plane-wave
Riemann problem and prove that the solution is unstable under
perturbations. Viscosity and heat conduction are taken into account
in section 3. In section 4 the relation to known nonphysical
numerical results (carbuncles) in Godunov´s and Roe´s method is
discussed and an explanation for carbuncle instabilities proposed,
which is also consistent with Majda´s shock stability analysis
[MA1983]. The last section contains a short summary and some
conclusions.
\\
\section{An Instability in the plane-wave \\ Riemann Problem}
In this section we analyze the stability of the plane-wave Riemann
problem. Since the Euler and Navier-Stokes equations are rotational
symmetric, we can assume without loss of generality that the plane
wave moves in the x-direction and that the initial states are
separated by the line $x(y)=0$.
Assume the solution of the plane-wave Riemann is given by a single
discontinuity, then the initial states $\mathbf{u}_l$ and
$\mathbf{u}_r$ of the Riemann problem satisfy the Rankine-Hugoniot
jump condition [CF1948].
\begin{align}
    \dot{s} [\mathbf{\mathbf{u}}_r - \mathbf{\mathbf{u}}_l]
    = f(\mathbf{\mathbf{u}}_r) - f(\mathbf{\mathbf{u}}_l) \label{eq:rankine}
\end{align}
where $\dot{s}$ is the shock speed.
The normal velocity of the discontinuity can be computed from
$\eqref{eq:rankine}$.
\newline
Let us assume that the states $\mathbf{u}_l$ and $\mathbf{u}_r$ are
connected by a near stationary 1-shock wave with normal shock speed
$\dot{s}^1=\dot{s}^1(\mathbf{u}_l,\mathbf{u}_r) < 0$. Since the
shock is near stationary, we have $-\varepsilon^s < \dot{s}^1$, for
a small positive number $\varepsilon^s$. Let
${\mathbf{\tilde{u}}}_r$ be a small perturbation of the right state.
The perturbed plane wave $\tilde{\mathbf{u}}$ has the initial data
\begin{align}
    \tilde{\mathbf{u}}(x,y,0)=
    \begin{cases}
        \mathbf{u}_l & \text{for  } x < 0 \label{eq:per-rie}\\
        \mathbf{\tilde{u}}_r & \text{for  } 0 < x
    \end{cases}
\end{align}
and consists of a 1-shock, a contact discontinuity and a 3-wave,
which is a weak rarefaction or weak shock wave.\newline In the
following we neglect, without loss of generality, the 3-wave such
that $\mathbf{\tilde{u}}_{mr}=\mathbf{\tilde{u}}_r$. A smooth
function F exists with
\begin{align}
    \mathbf{\tilde{\mathbf{u}}}_r= \mathbf{u}_l +
    F(\varepsilon_1,\varepsilon_2,\varepsilon_3;\mathbf{u}_l) \label{eq:per-rie-F}
\end{align}
where $\varepsilon_2$ represents a parameter for the strength of the
contact discontinuity and $\varepsilon_3$ is a parameter for the
strength of a 3-wave; see [SM1983; Chapter 17]. $\varepsilon_1$
represents a parameter for the strength of a 1-shock.
We choose $\varepsilon_1$ such that
\begin{align}
    \mathbf{u}_r = \mathbf{\mathbf{u}}_l + F(\varepsilon_1,0,0;\mathbf{u}_l)
\end{align}
For a perturbation of $\varepsilon_2$ with $\varepsilon_3 = 0$ the
pressure p and the x-component of the velocity are constant behind
the 1-shock. We obtain form the Lax shock conditions
$\dot{s}^1_{i+1/2,j} < \lambda_2(\mathbf{u}_r)$ where
$\lambda_2(\mathbf{u}_r)$ is the second eigenvalue of the Jacobian
$A(\mathbf{u}_r)$.
For a near stationary strong shock wave we can assume that $0 <
\lambda_2(\mathbf{u}_r)$. Since the parameter $\varepsilon_2$
represents a contact discontinuity, the discontinuity speed
$\dot{s}^2$ is equal to the characteristic speed
$\lambda_2(\mathbf{u}_r)$, which is constant across a contact
discontinuity; i.e.
$\dot{s}^2=\lambda_2(\mathbf{u}_r)=\lambda_2(\mathbf{\tilde{u}}_{r})$.
Where $\mathbf{\tilde{u}}_{r}$ is given by
\begin{align}
    \mathbf{\tilde{u}}_{r}=\mathbf{u}_l + F(\varepsilon_1,\varepsilon_2,0;\mathbf{u}_l) \label{eq:rp-def-1}
\end{align}
The y-component of the velocity $v_l$ and $v_r$ enters the solution
of the plane-wave Riemann problem essentially as a parameter. We can
first solve the one dimensional Riemann problem ignoring the
momentum equation for $v$ and then introduce a jump in $v$ at the
contact discontinuity to obtain the full plane-wave solution.
However, $v_l$ and $v_r$ enter the one dimensional Riemann problem
through the pressure as a function of the total energy, density and
velocity. We assume that $v_r=O(\varepsilon)$ and $v_l=
O(\varepsilon)$. For a small perturbation $O(\varepsilon)$ the
change in $\eqref{eq:eqof}$ is of order $\varepsilon^2$. Therefore
we can assume that the wave structure for the perturbed plane-wave
Riemann problem
\begin{align}
    \tilde{\mathbf{u}}(x,y,0)= \label{eq:per-rie1}
    \begin{cases}
        \mathbf{u}_l^\varepsilon (0-,y)  & \text{for  } x < 0 \\
        \mathbf{\tilde{u}}_r^\varepsilon (0+,y) & \text{for  } 0 < x
    \end{cases}
\end{align}
persists up to second order in $\varepsilon$, where
$\mathbf{u}^\varepsilon_l (x,y) =
(\rho_l,\rho_lu_l,\rho_lv^{\varepsilon}(x,y),E_l)^T$ and
$\mathbf{u}_r^\varepsilon (x,y) =
(\tilde{\rho}_r,\tilde{\rho}_r{u}_r,\tilde{\rho_r}v^{\varepsilon}(x,y),\tilde{E}_r)^T$.
$\\$
\newline If $v^{\varepsilon}(x,y)$ is
discontinuous at $x=0$ the perturbation introduces a jump in the
y-component of the velocity at the contact discontinuity. If
$v_r\neq v_l$ the perturbed plane-wave Riemann problem contains a
tangential or shear instability. In [LL1959;§81] it is proved that
such a tangential instability in an incompressible non viscous flow
is absolutely unstable and may lead to a turbulent flow, and it is
further mentioned that these instabilities also exists in
compressible flows. However, we assume in the following that the
perturbation $v^{\varepsilon}(x,y)$ is a smooth function.
\newline
We can consider the parameter $\varepsilon_1$, $\varepsilon_2$ in
$\eqref{eq:rp-def-1}$ as functions of $y$. We assume that the
variation of $\varepsilon_1$ in $y$ is small enough such that
$\dot{s}^1_{i+1/2,j} < 0 < u_r$ still holds. This variation can be
defined independently of the y-component of the velocity and we can
assume that $v$ is not affected through this perturbation.
The flow for this perturbed plane-wave Riemann problem is defined
through the two-dimensional Euler equation $\eqref{eq:Euler-diff}$.
A change of the parameter $\varepsilon_2$ in $\eqref{eq:rp-def-1}$
affects only the contact discontinuity. Since the pressure and the
x-component of the velocity are constant across a contact
discontinuity, we have behind the 1-shock
$p(x,y,t)=p_r=\tilde{p}_r$, $u(x,y,t)=u_r=\tilde{u}_r$. Thus the
Euler equations reduce behind the 1-shock to:
\begin{subequations}
\begin{align}
    &\frac{\partial}{\partial t}\rho + u_r \frac{\partial}{\partial x}\rho + v \frac{\partial}{\partial y}\rho + \rho \frac{\partial}{\partial y}v = 0 \label{eq:shear-0-a} \\
    &\frac{\partial}{\partial t}v + u_r\frac{\partial}{\partial x}v + v\frac{\partial}{\partial y}v= 0
     \label{eq:shear-0-b} \\
    & \frac{\partial}{\partial t}E + u_r \frac{\partial}{\partial x}E + v\frac{\partial}{\partial
    y}E + E\frac{\partial}{\partial y}v = 0 \label{eq:shear-0-c}
\end{align} \label{eq:shear-0-flow}
\end{subequations}
\newline
Using $\eqref{eq:E=rhop}$ to rewrite the total energy in
$\eqref{eq:Euler-diff}$ as a function of the density, pressure and
velocities and using the other conservative equations in
$\eqref{eq:Euler-diff}$, the energy equation may be rewritten as a
pressure equation
\begin{align}
    &\frac{\partial}{\partial t}p + u_r \frac{\partial}{\partial x}p + v \frac{\partial}{\partial y}p
    + \gamma p(\frac{\partial}{\partial x}u + \frac{\partial}{\partial y}v) = 0 \label{eq:shear-p}
\end{align}
which for a constant pressure simply reduce to a divergence
condition
\begin{align}
    \frac{\partial}{\partial x}u + \frac{\partial}{\partial y}v=0
\end{align}
Since the x-component of the velocity is constant behind the 1-Shock
$\eqref{eq:shear-0-flow}$ reduces to:
\begin{subequations}
\begin{align}
    &\frac{\partial}{\partial t}\rho + u_r \frac{\partial}{\partial x}\rho + v \frac{\partial}{\partial y}\rho = 0 \label{eq:shear-f-a} \\
    &\frac{\partial}{\partial t}v + u_r\frac{\partial}{\partial x}v = 0
     \label{eq:shear-f-b} \\
    & \frac{\partial}{\partial t}E + u_r \frac{\partial}{\partial x}E + v\frac{\partial}{\partial
    y}E = 0 \label{eq:shear-f-c}
\end{align} \label{eq:shear-flow}
\end{subequations}
\newline
where the last equation follows from the equation of state
$\eqref{eq:eqof}$ and from the first two. Therefore the perturbed
solution satisfies behind the 1-shock the equations
$\eqref{eq:shear-flow}$ with the initial data
\begin{subequations}
\begin{align}
    &\mathbf{u}(x,y,0+)= \label{eq:rie-shear-int-a}
    \begin{cases}
        \mathbf{u}^{\epsilon}(x,y) & \text{for  } x < 0 \\
        \mathbf{\tilde{u}}^{\epsilon}(x,y) & \text{for  } 0 < x
    \end{cases} \\
    &\text{with} \nonumber \\
    &\mathbf{u}^{\epsilon}(x,y)=(\rho_0(x,y),\rho_0(x,y) u_r,\rho_0(x,y) v^{\varepsilon}_0(x,y),
    E_0(x,y))^T \label{eq:rie-shear-int-rhol}    \text{  \;  } \\
    &\text{and} \nonumber \\
    &\mathbf{\tilde{u}}^{\epsilon}(x,y)=(\tilde{\rho}_0(x,y),\tilde{\rho}_0(x,y) u_r,\tilde{\rho}_0(x,y) v^{\varepsilon}_0(x,y), E_0(x,y))^T
    \text{   \;  } \label{eq:rie-shear-int-rhor} \\
    &\text{and} \nonumber \\
    &E_0(x,y)= \label{eq:rie-shear-int-E}
    \begin{cases}
        p_r/(\gamma-1)+\frac{1}{2}\rho_0(x,y)(u^2_r+v^{\varepsilon}_0(x,y)^2) & \text{for  } x < 0 \\
        p_r/(\gamma-1)+\frac{1}{2}\tilde{\rho}_0(x,y)(u^2_r+v^{\varepsilon}_0(x,y)^2) & \text{for  } 0 < x
    \end{cases}
\end{align} \label{eq:rie-shear-int}
\end{subequations}
where $v_0^{\varepsilon}(x,y)$, $\rho_0(x,y)$ and
$\tilde{\rho}_0(x,y)$ are smooth perturbed initial data.
\\
\\
$\mathbf{Proposition\;1}$: a) If the initial tangential velocity
$v_0^{\varepsilon}$ does not depend on $y$ then the constant
pressure solution of the initial value problem
$\eqref{eq:shear-flow}$, $\eqref{eq:rie-shear-int}$ is given by:
\begin{subequations}
\begin{align}
    &\rho(x,y)= \label{eq:prop1-rho}
    \begin{cases}
        \rho_0(x - u_r t,y - v^{\varepsilon} t) & \text{for  } x < u_r t \\
        \tilde{\rho}_0(x - u_r t,y - v^{\varepsilon} t) & \text{for  } 0 < u_r t
    \end{cases} \\
    &\text{and} \nonumber \\
    &v^{\varepsilon} = v_0^{\varepsilon}(x - u_r t)
    \label{eq:prop1-v} \\
    &\text{and} \nonumber \\
    &E = p_r/(\gamma-1)+\frac{1}{2}\rho(x,y)(u^2_r+(v^{\varepsilon})^2) \label{eq:prop1-E}
\end{align} \label{eq:shea-flow-label}
\end{subequations}
b) If the initial tangential velocity $v_0^{\varepsilon}$ depends on
$y$ a solution does not exist and the pressure is not constant.
\newline
\newline
Proof: The solution of the plane-wave Riemann problem has behind the
1-shock, i.e. for $x > s^1_{i+1/2,j}t$, a constant pressure $p_r$
and a constant x-component of the velocity $u_r$. Therefore the
Euler equations reduces to $\eqref{eq:shear-flow}$ and the solution
is completely defined through the density $\rho$, the tangential
velocity $v$, the constant pressure and the constant x-component of
the velocity.
We have for $x < u_r t$
\begin{align}
    &\frac{\partial}{\partial t}\rho + u_r \frac{\partial}{\partial
    x}\rho + v^{\varepsilon}\frac{\partial}{\partial y}\rho \nonumber \\
    &=[-u_r + u_r]\frac{\partial}{\partial x}\rho_0
    +[-(v^{\varepsilon} + t \frac{\partial}{\partial t}v^{\varepsilon}) - u_r t \frac{\partial}{\partial
    x}v^{\varepsilon} + v^{\varepsilon}( 1 - t\frac{\partial}{\partial y}v^{\varepsilon})]\frac{\partial}{\partial y}\rho_0
    \nonumber \\
    &=-[\frac{\partial}{\partial t}v^{\varepsilon} + u_r \frac{\partial}{\partial
    x}v^{\varepsilon} + \frac{\partial}{\partial y}v^{\varepsilon}]t\frac{\partial}{\partial y}\rho_0 \label{eq:shea-flow-proof-1}
\end{align}
and similar for $x > u_r t$
\begin{align}
    &\frac{\partial}{\partial t}\tilde{\rho} + u_r \frac{\partial}{\partial
    x}\tilde{\rho} + v^{\varepsilon}\frac{\partial}{\partial y}\tilde{\rho} \nonumber \\
    &=[-u_r + u_r]\frac{\partial}{\partial x}\tilde{\rho}_0
    +[-(v^{\varepsilon} + t \frac{\partial}{\partial t}v^{\varepsilon}) - u_r t \frac{\partial}{\partial
    x}v^{\varepsilon} + v^{\varepsilon}( 1 -  v^{\varepsilon} t\frac{\partial}{\partial y}v^{\varepsilon})]\frac{\partial}{\partial y}\tilde{\rho}_0
    \nonumber \\
    &=-[\frac{\partial}{\partial t}v^{\varepsilon} + u_r \frac{\partial}{\partial
    x}v^{\varepsilon} +  v^{\varepsilon} \frac{\partial}{\partial y}v^{\varepsilon}]t\frac{\partial}{\partial y}\tilde{\rho}_0 \label{eq:shea-flow-proof-2}
\end{align}
Furthermore
\begin{align}
    &\frac{\partial}{\partial t}v^{\varepsilon} + u_r \frac{\partial}{\partial x}v^{\varepsilon}
    =[-u_r + u_r]\frac{\partial}{\partial x}v^{\varepsilon}_0 = 0
\end{align}
If follows form $\eqref{eq:shea-flow-proof-1}$,
$\eqref{eq:shea-flow-proof-2}$ for an initial tangential velocity
$v^{\varepsilon}_0=v^{\varepsilon}_0(x)$ that $\eqref{eq:prop1-rho}$
and $\eqref{eq:prop1-v}$ are solutions of the initial value problem
$\eqref{eq:shear-flow}$, $\eqref{eq:rie-shear-int}$.
\newline
We obtain furthermore for the total energy E for $x \neq u_r t$:
\begin{align}
    & \frac{\partial}{\partial t}E + u_r \frac{\partial}{\partial x}E + v^{\varepsilon}\frac{\partial}{\partial y}E
    \nonumber \\
    &=\frac{1}{2} (u_r^2 +(v^{\epsilon})^2)[\frac{\partial}{\partial t}\rho + u_r \frac{\partial}{\partial x}\rho + v^{\varepsilon} \frac{\partial}{\partial
    y}\rho]\nonumber \\
    &+\rho v^{\varepsilon}[\frac{\partial}{\partial t}v^{\varepsilon} + u_r\frac{\partial}{\partial x}v^{\varepsilon}]=0 \nonumber
\end{align}
Therefore $\eqref{eq:shea-flow-label}$ defines a solution of the
Euler equation $\eqref{eq:Euler-diff}$ in smooth parts of the flow.
\newline
The jump condition for a plane-wave moving in the x-direction
reduces for a constant pressure and constant x-component of the
velocity to
\begin{align}
    &\dot{s}(t)[\rho(x+,y,t)-\rho(x-,y,t)]=u_r[\rho(x+,y,t)-\rho(x-,y,t)] \nonumber \\
    &\dot{s}(t)u_r[\rho(x+,y,t)-\rho(x-,y,t)]=(u_r)^2[\rho(x+,y,t)-\rho_r(x-,y,t)] \nonumber \\
    &\dot{s}(t)[\rho(x+,y,t)v^{\varepsilon}(x+,t)-\rho(x-,y,t)v^{\varepsilon}(x-,t)]= \nonumber \\
    &u_r[\rho(x+,y,t)v^{\varepsilon}(x+,t)-\rho(x-,y,t)v^{\varepsilon}(x-,t)] \nonumber \\
    &\dot{s}(t)[E(x+,y,t)-E(x-,y,t)]=u_r[E(x+,y,t)-E(x-,y,t)] \nonumber
\end{align}
where $\dot{s}$ is the speed of the discontinuity in the
x-direction. Therefore we see that the jump conditions are satisfied
along the curve $x=s(t)=u_rt$, and $\eqref{eq:shea-flow-label}$ is a
weak solution of $\eqref{eq:Euler-diff}$.\\
Since the solution of initial value problem $\eqref{eq:shear-flow}$,
$\eqref{eq:rie-shear-int}$ is unique for smooth velocities $u_r$ and
$v^{\varepsilon}$, part a) of the proposition follows.
\\
If the initial data $v^{\varepsilon}_0$ depend on $y$ the unique
solution of $\eqref{eq:shear-f-b}$ is for every $y$ given by
$\eqref{eq:prop1-v}$. In this case we obtain for the total energy E
for $x \neq u_r t$ and a constant pressure:
\begin{align}
    & \frac{\partial}{\partial t}E + u_r \frac{\partial}{\partial x}E + v^{\varepsilon}\frac{\partial}{\partial y}E
    \nonumber \\
    &=\frac{1}{2} (u_r^2 +(v^{\epsilon})^2)[\frac{\partial}{\partial t}\rho + u_r \frac{\partial}{\partial x}\rho + v^{\varepsilon} \frac{\partial}{\partial
    y}\rho]\nonumber \\
    &+\rho v^{\varepsilon}[\frac{\partial}{\partial t}v^{\varepsilon} + u_r\frac{\partial}{\partial x}v^{\varepsilon}
    + v^{\varepsilon} \frac{\partial}{\partial y}v^{\varepsilon}]
    \nonumber \\
    &=\rho v^{\varepsilon} v^{\varepsilon} \frac{\partial}{\partial y}v^{\varepsilon}
\end{align}
Thus the energy equation can only be satisfied if $v^{\varepsilon}$
does not depend on $y$. This must also hold for $t=0$. If $v$
depends on $y$ we obtain from $\eqref{eq:shear-p}$
\begin{align}
    &\frac{\partial}{\partial t}p + u_r \frac{\partial}{\partial x}p + v \frac{\partial}{\partial y}p
    = - \gamma p \frac{\partial}{\partial y}v  \label{eq:prop1-shear-p}
\end{align}
which cannot be satisfied for a constant pressure solution. This
proves part b) and completes the proof. $\Box$
\\
\newline In a constant pressure flow the equations for the conservation of energy follows from the
momentum and the continuity equations. The Euler equations
$\eqref{eq:Euler-diff}$ therefore simplify to
\begin{subequations}
\begin{align}
    &\frac{\partial}{\partial t}\rho + u \frac{\partial}{\partial
    x}\rho  + v \frac{\partial}{\partial y}\rho  = 0 \label{eq:pconst-3a} \\
   & \frac{\partial}{\partial t} u + u \frac{\partial}{\partial x} u + v \frac{\partial}{\partial y} u = 0 \label{eq:pconst-3b} \\
    &\frac{\partial}{\partial t} v + u \frac{\partial}{\partial x} v + v \frac{\partial}{\partial y} v = 0 \label{eq:pconst-3c}
\end{align} \label{eq:pconst}
\end{subequations}
\newline
If the velocity field is given by $u=u_r$ and $v=v(x,t)$ these
equations reduce to the equations $\eqref{eq:shear-flow}$ for a
perturbed flow behind a shock. This system of partial differential
equations is unstable under perturbations. Let $\bar{\rho}$,
$\bar{u}$ and $\bar{v}$ be a constant pressure solution of
$\eqref{eq:pconst}$ and $\tilde{\rho}$, $\tilde{u}$, $\tilde{v}$
infinitesimal perturbations. It is assumed that
\begin{align}
    \rho= \bar{\rho} + \tilde{\rho} \;\;\;\;u=\bar{u} + \tilde{u} \;\;\;\;
    v= \bar{v} + \tilde{v} \nonumber
\end{align}
is also a smooth constant pressure solution of $\eqref{eq:pconst}$
$\footnote{The results remain valid if pressure perturbations
$p=\bar{p} + \tilde{p}$  are included}.$ We obtain by neglecting
quadratic perturbation terms
\begin{align}
 & \frac{\partial}{\partial t}\bar{\rho} + \bar{u} \frac{\partial}{\partial
 x}\bar{\rho} +  \bar{v} \frac{\partial}{\partial y}\bar{\rho}
 + \frac{\partial}{\partial t}\tilde{\rho} + \bar{u} \frac{\partial}{\partial x}\tilde{\rho}
 +  \bar{v} \frac{\partial}{\partial y}\tilde{\rho} + \tilde{u} \frac{\partial}{\partial
 x}\bar{\rho} +  \tilde{v} \frac{\partial}{\partial y}\bar{\rho} = 0
 \nonumber \\
 & \frac{\partial}{\partial t}\bar{u} + \bar{u} \frac{\partial}{\partial
 x}\bar{u} +  \bar{v} \frac{\partial}{\partial y}\bar{u}
 + \frac{\partial}{\partial t}\tilde{u} + \bar{u} \frac{\partial}{\partial x}\tilde{u}
 +  \bar{v} \frac{\partial}{\partial y}\tilde{u} + \tilde{u} \frac{\partial}{\partial
 x}\bar{u} +  \tilde{v} \frac{\partial}{\partial y}\bar{u} = 0
 \nonumber \\
 & \frac{\partial}{\partial t}\bar{v} + \bar{u} \frac{\partial}{\partial
 x}\bar{v} +  \bar{v} \frac{\partial}{\partial y}\bar{v}
 + \frac{\partial}{\partial t}\tilde{v} + \bar{u} \frac{\partial}{\partial x}\tilde{v}
 +  \bar{v} \frac{\partial}{\partial y}\tilde{v} + \tilde{u} \frac{\partial}{\partial
 x}\bar{v} +  \tilde{v} \frac{\partial}{\partial y}\bar{v} = 0
 \nonumber
\end{align}
which can be rewritten as a inhomogeneous hyperbolic system
\begin{subequations}
\begin{align}
 & \frac{\partial}{\partial t}\tilde{\rho} + \bar{u} \frac{\partial}{\partial x}\tilde{\rho}
 +  \bar{v} \frac{\partial}{\partial y}\tilde{\rho} = - \tilde{u} \frac{\partial}{\partial
 x}\bar{\rho} -  \tilde{v} \frac{\partial}{\partial y}\bar{\rho}
 \label{eq:pw-rho2} \\
 & \frac{\partial}{\partial t}\tilde{u} + \bar{u} \frac{\partial}{\partial x}\tilde{u}
 +  \bar{v} \frac{\partial}{\partial y}\tilde{u} = - \tilde{u} \frac{\partial}{\partial
 x}\bar{u} -  \tilde{v} \frac{\partial}{\partial y}\bar{u}
 \label{eq:pw-nx2} \\
 & \frac{\partial}{\partial t}\tilde{v} + \bar{u} \frac{\partial}{\partial x}\tilde{v}
 +  \bar{v} \frac{\partial}{\partial y}\tilde{v} = - \tilde{u} \frac{\partial}{\partial
 x}\bar{v} -  \tilde{v} \frac{\partial}{\partial y}\bar{v}
 \label{eq:pw-ny2}
\end{align} \label{eq:pw-main}
\end{subequations}
The solution of the homogeneous system $\eqref{eq:pw-main}$ is
\begin{subequations}
\begin{align}
    \tilde{\rho}^h=&\tilde{\rho}^0(x-\bar{u}t, y-\bar{v}t) \label{eq:pw-rho3} \\
    \tilde{u}^h=&\tilde{u}^0(x-\bar{u}t, y-\bar{v}t) \label{eq:pw-nx3} \\
    \tilde{v}^h=&\tilde{v}^0(x-\bar{u}t, y-\bar{v}t) \label{eq:pw-ny3}
\end{align}
\end{subequations} where $\tilde{\rho}^0$, $\tilde{u}^0$ and
$\tilde{v}^0$ are some arbitrary initial divergence free
perturbations. Let
\begin{subequations}
\begin{align}
    S^{\rho}(x,y,t)= - \tilde{u}^h \frac{\partial}{\partial x}\bar{\rho}
                - \tilde{v}^h \frac{\partial}{\partial y}\bar{\rho}  \label{eq:pw-source-rho} \\
    S^x(x,y,t)= - \tilde{u}^h \frac{\partial}{\partial x}\bar{u}
                - \tilde{v}^h \frac{\partial}{\partial y}\bar{u}  \label{eq:pw-source-x} \\
    S^y(x,y,t)= - \tilde{u}^h \frac{\partial}{\partial x}\bar{v}
                - \tilde{v}^h \frac{\partial}{\partial y}\bar{v}  \label{eq:pw-source-y}
\end{align} \label{eq:pw-source-main}
\end{subequations} \newline
By formally applying Duhamel´s principle we obtain for the
inhomogeneous equation $\eqref{eq:pw-main}$ with the source term
$\eqref{eq:pw-source-main}$
\begin{subequations}
\begin{align}
    \tilde{\rho}   &= \tilde{\rho}^h + \int_0^t S^{\rho}(x - \bar{u}(t-\tau), y - \bar{v}(t-\tau),\tau) d\tau \nonumber \\
                &=\tilde{\rho}^h - t(\tilde{u}^h \frac{\partial}{\partial x}\bar{\rho} + \tilde{v}^h \frac{\partial}{\partial y}\bar{\rho}) \label{eq:pw-rho6} \\
    \tilde{u}   &= \tilde{u}^h + \int_0^t S^x(x - \bar{u}(t-\tau), y - \bar{v}(t-\tau),\tau) d\tau \nonumber \\
                &=\tilde{u}^h - t(\tilde{u}^h \frac{\partial}{\partial x}\bar{u} + \tilde{v}^h \frac{\partial}{\partial y}\bar{u}) \label{eq:pw-nx6} \\
    \tilde{v}   &= \tilde{v}^h + \int_0^t S^y(x - \bar{u}(t-\tau), y - \bar{v}(t-\tau),\tau) d\tau \nonumber \\
                &= \tilde{v}^h - t(\tilde{u}^h \frac{\partial}{\partial x}\bar{v} + \tilde{v}^h \frac{\partial}{\partial y}\bar{v}) \label{eq:pw-ny7}
\end{align}\label{eq:pw-nyx}
\end{subequations}
\newline
where $\bar{\rho}=\bar{\rho}(x,y,t)$ , $\bar{u}=\bar{u}(x,y,t)$ and
$\bar{v}=\bar{v}(x,y,t)$. Inserting $\tilde{\rho}$, $\tilde{u}$ and
$\tilde{v}$ into $\eqref{eq:pconst}$ we can verify that
$\eqref{eq:pw-nyx}$ defines a solution of the inhomogeneous
hyperbolic system $\eqref{eq:pw-main}$. We obtained, that for a non
constant background flow $\bar{\rho}$, $\bar{u}$ and $\bar{v}$,
perturbation $\tilde{\rho}$, $\tilde{u}$ and $\tilde{v}$ grow linear
in time.
\newline
The inhomogeneous hyperbolic system $\eqref{eq:pw-main}$ can be
rewritten in vector form as
\begin{align}
   & \frac{\partial}{\partial t} \mathbf{\tilde{w}}(x,y,t)
   + \mathbf{\bar{A}}\frac{\partial}{\partial x}\mathbf{\tilde{w}}(x,y,t)
   + \mathbf{\bar{B}} \frac{\partial}{\partial y} \mathbf{\tilde{w}}(x,y,t)
   + \mathbf{\bar{C}}\mathbf{\tilde{w}}(x,y,t)= 0\label{eq:Euler-pert-110}
\end{align}
where the vector of primitive variables is
$\mathbf{\tilde{w}}=(\tilde{\rho},\tilde{u},\tilde{v})^T$ and the 3
x 3 matrices $\bar{A}(x,y,t)$, $\bar{B}(x,y,t)$ and $\bar{C}(x,y,t)$
are defined by
\begin{align}
\mathbf{\bar{A}}=\begin{pmatrix}
 \bar{u} & 0 & 0\\
  0  & \bar{u} & 0  \\
  0 & 0 & \bar{u}
\end{pmatrix}
\;\;\;\;\;\mathbf{\bar{B}}=\begin{pmatrix}
 \bar{v} & 0 & 0 \\
  0 &  \bar{v} & 0 \\
  0 & 0 & \bar{v}
\end{pmatrix}
\;\;\;\;\;-\mathbf{\bar{C}}=\begin{pmatrix}
 0 & \frac{\partial}{\partial x}\bar{\rho} & \frac{\partial}{\partial y}\bar{\rho}   \\
 0 & \frac{\partial}{\partial x}\bar{u} & \frac{\partial}{\partial y}\bar{u} \\
 0 & \frac{\partial}{\partial x}\bar{v} & \frac{\partial}{\partial
 y}\bar{v}
\end{pmatrix} \nonumber
\end{align} $\label{eq:pw-main-matrix}$
\newline
For a smooth background flow $\eqref{eq:Euler-pert-110}$ is a
symmetric hyperbolic system, with matrices $\bar{A}(x,y,t)$,
$\bar{B}(x,y,t)$ and $\bar{C}(x,y,t)$ depending smoothly on $x$, $y$
and $t$. For any perturbation $\mathbf{\tilde{w}}$ with compact
support in $(x,y)$ a energy inequality of the form
\begin{align}
    \parallel \mathbf{\tilde{w}}(t) \parallel \leq \exp(M \mid t \mid)\parallel \mathbf{\tilde{w}}(0)
    \parallel \label{eq:pw-energy-inequality}
\end{align}
holds, where the constant $M$ depends only on the magnitude of the
symmetric part of $\bar{C}$ and of the first derivative of $\bar{A}$
with respect to $x$ resp. $\bar{B}$ with respect to $y$; see
[LAX2006; Section 4.3]. This result can also be derived directly
from $\eqref{eq:pw-nyx}$. Now consider the solution
$\eqref{eq:pw-nyx}$ for a case where the initial background density
$\bar{\rho}$ is not differentiable at $x=0$. For a constant
x-component of the velocity $\bar{u}=u_r$ and a initial tangential
velocity $\bar{v}_0=\bar{v}_0(x)$, which is enforced through a
constant pressure, the solution is given by
$\eqref{eq:shea-flow-label}$. In this case the definition of the
perturbed solution $\eqref{eq:pw-nyx}$ remains valid, if we
interpret the partial derivative of $\bar{\rho}$ with respect to $x$
at the contact discontinuity as a delta function; i.e.
\begin{align}
    \frac{\partial}{\partial x}\bar{\rho}(u_r t, y, t)=(\bar{\rho}_r
    - \bar{\rho}_l)\delta(x - u_r t)
\end{align}
where $\bar{\rho}_l$ and $\bar{\rho}_r$ denote the left and right
limits at the contact discontinuity. For an initial oscillatory
velocity perturbation $\tilde{u}^h$, large density perturbations are
enforced. In this case the solution grows with $t$ in a manner
unlike what we would expect from a hyperbolic equation and no energy
inequality of the form $\eqref{eq:pw-energy-inequality}$ can be
derived. We denote such a flow in the following as unstable under
perturbations.
\newline Consider for example a vanishing background tangential velocity
 $\bar{v}$, a vanishing perturbed tangential velocity $\tilde{v}$ and a
piecewise constant background density $\bar{\rho}(x,t)$ which is
discontinuous at $x=\bar{u}t$, then we obtain from
 $\eqref{eq:pw-rho2}$ for $\bar{u}=u_r$,
\begin{align}
 & \frac{\partial}{\partial t}\tilde{\rho} + u_r \frac{\partial}{\partial x}\tilde{\rho} = - \tilde{u} (\bar{\rho}_r
    - \bar{\rho}_l)\delta(x - u_r t) \label{eq:epw-simple}
\end{align}
and from $\eqref{eq:pw-rho6}$ for the solution
\begin{align}
    \tilde{\rho}=\tilde{\rho}^h - t \tilde{u}^h (\bar{\rho}_r - \bar{\rho}_l)\delta(x - u_r
    t) \label{eq:epw-nyx-sol}
\end{align}
For an initial oscillatory velocity perturbation $\tilde{u}^h$ large
density perturbations are enforced.
\newline For a constant normal
velocity $\bar{u}=u_r$ and a piecewise constant background
tangential velocity
 $\bar{v}$ which depends only on x and is discontinuous at $x=u_r
 t$, equation $\eqref{eq:pw-ny2}$ becomes
 \begin{align}
 & \frac{\partial}{\partial t}\tilde{v} + u_r \frac{\partial}{\partial x}\tilde{v}
   = - \tilde{u} (\bar{v}_r - \bar{v}_l)\delta(x - u_r t) \label{eq:epw-simple-1}
 \end{align}
 with a weak solution given by $\eqref{eq:pw-ny7}$; i.e.
 \begin{align}
    \tilde{v}=\tilde{v}^h - t \tilde{u}^h (\bar{v}_r - \bar{v}_l)\delta(x - u_r
    t) \label{eq:epw-nyx-sol1}
\end{align}
For an initial oscillatory velocity perturbation $\tilde{u}^h$,
large tangential velocity perturbations are enforced. Linear
advection equation of the form $\eqref{eq:epw-simple}$,
$\eqref{eq:epw-simple-1}$ with a singular source term are discussed
by LeVeque in [LEV2002; Section 16.3.1].
\newline
\newline
We showed that the flow behind a perturbed plane stationary shock
line is a constant pressure flow governed by $\eqref{eq:pconst}$
with $u=u_r$ and $v=v(x,t)$. Generally, the flow between the
nonlinear 1-wave and 3-wave in the Riemann problem is a constant
pressure flow region. Since a discontinuity in a constant pressure
plane wave flow region is a contact discontinuity, we obtained:
\newline
\newline
\emph{If the solution of the plane-wave Riemann problem contains a
contact discontinuity, then the solution is unstable under
perturbations.}
\\
\\
Remark: The same perturbation analysis can be applied to an
homentrop (constant entropy), nearly constant density flow, with a
not necessary constant pressure and divergence. In this case we
obtain for the perturbation equations
\begin{subequations}
\begin{align}
 & \frac{\partial}{\partial t}\tilde{\rho} + \bar{u}
\frac{\partial}{\partial x}\tilde{\rho}
 +  \bar{v} \frac{\partial}{\partial y}\tilde{\rho} + \bar{\rho}_0
(\frac{\partial}{\partial x}\tilde{u}
 +  \frac{\partial}{\partial y}\tilde{v}) = 0
 \label{eq:apw-rho6m} \\
 & \frac{\partial}{\partial t}\tilde{u} + \bar{u} \frac{\partial}{\partial x}\tilde{u}
 +  \bar{v} \frac{\partial}{\partial y}\tilde{u} = - \tilde{u} \frac{\partial}{\partial
 x}\bar{u} -  \tilde{v} \frac{\partial}{\partial y}\bar{u} - \frac{1}{\bar{\rho}_0} \frac{\partial}{\partial x}\tilde{p}
 \label{eq:apw-nx6m} \\
 & \frac{\partial}{\partial t}\tilde{v} + \bar{u} \frac{\partial}{\partial x}\tilde{v}
 +  \bar{v} \frac{\partial}{\partial y}\tilde{v} = - \tilde{u} \frac{\partial}{\partial
 x}\bar{v} -  \tilde{v} \frac{\partial}{\partial y}\bar{v} - \frac{1}{\bar{\rho}_0} \frac{\partial}{\partial y}\tilde{p}
 \label{eq:apw-ny6m} \\
  & \frac{\partial}{\partial t}\tilde{p} + \bar{u} \frac{\partial}{\partial x}\tilde{p}
 +  \bar{v} \frac{\partial}{\partial y}\tilde{p} + \gamma\bar{p}
(\frac{\partial}{\partial x}\tilde{u}
 +  \frac{\partial}{\partial y}\tilde{v}) =  0
 \label{eq:apw-p6m}
\end{align} \label{eq:apw-main6m}
\end{subequations}
\newline where $\bar{\rho}_0$ is the constant background density and
$(\bar{u},\bar{v})$ a divergence free background velocity field. For
a constant background flow $\bar{\rho}_0, \bar{u}_0,\bar{v}_0,
\bar{p}_0$, homentrop density perturbations are defined through
\begin{align}
    \tilde{\rho}=\tilde{\rho}_A(x-\bar{u}_0t,y-\bar{v}_0t,t) \label{eq:isenf-per-rho}
\end{align}
where $\tilde{\rho}_A$ satisfies a wave equation
\begin{align}
    \frac{\partial^2}{\partial \tau^2}\tilde{\rho}_A - \bar{c}_0^2
    ( \frac{\partial^2}{\partial x^2}\tilde{\rho}_A + \frac{\partial^2}{\partial
    y^2}\tilde{\rho}_A) = 0 \label{eq:wave-eq}
\end{align}
where $\bar{c}_0$ denotes the sound speed of the background flow and
$\tau$ refers to the differentiation with respect to the third
argument in $\tilde{\rho}_A$. Since the solutions of the wave
equation are stable, density perturbations are stable.
\newline
A initial (acoustic) perturbation $\tilde{\rho}_A$ propagates with
sound speed relative to the background velocity into the flow. If
generated at a shock, these acoustic perturbation propagate with
sound speed relative to adjacent (constant) velocities, into the
flow on both sides of the discontinuity. In contrast to acoustic
perturbations, density perturbations in a constant pressure region
can only exists on one side of a shock line. Since in an ideal gas,
any density perturbation in a constant pressure region is equivalent
to a entropy perturbation, the latter disturbances are also denoted
as entropy disturbances.
\newline
In the analysis of the perturbed Riemann problem, the weak 3-wave
was neglected. Since the change of the entropy across a shock wave
is of third order in $\varepsilon_3$ in $\eqref{eq:per-rie-F}$, the
weak 3-wave can be regarded as a (discontinuous) homentrop
perturbation of a constant state. The density of these acoustic
perturbation is defined through $\eqref{eq:isenf-per-rho}$,
$\eqref{eq:wave-eq}$ and do not cause an unstable flow.
\newline
\newline
Note: To derive the Rankine-Hugoniot jump condition
$\eqref{eq:rankine}$ from the integral form $\eqref{eq:Euler-int}$
of the conservation law, it is assumed that the integral
\begin{align}
    \int_{-\Delta_x / 2 }^{\Delta_x / 2} \frac{\partial}{\partial t} \mathbf{u}(\xi,y,t^n) \:
    d\xi \label{eq:Euler-int100}
\end{align}
approaches zero for $\lim \Delta_x \rightarrow0$. This assumption
fails for the perturbed two dimensional plane-wave Riemann problem.
Consider a discontinuity line $\mathcal{S}$. Without restriction of
generality we may assume that for a sufficiently small portion, the
discontinuity line $\mathcal{S}$ is perpendicular to the x-axis and
we may assume that the flow is smooth in the y-direction.
Furthermore for a sufficiently small time interval the shock speed
$\dot{s}$ may be considered as constant. Let $x=s(t,y)=\dot{s}(y)t$
be the spatial-time discontinuity surface across which $\mathbf{u}$
has a jump. We obtain from the conservation law
$\eqref{eq:Euler-diff}$ for $\Delta_x > 0$
\begin{align}
    \frac{\partial}{\partial t} \int_{-\Delta_x/2}^{+\Delta_x/2}\mathbf{u}(x,y,t) \: dx
    =& - [\mathbf{f}(\mathbf{u}(\Delta_x/2,y,t)) -
    \mathbf{f}(\mathbf{u}(-\Delta_x/2,y,t))] \nonumber \\
    &-  \int_{-\Delta_x/2}^{+\Delta_x/2} \frac{\partial}{\partial y} \mathbf{g}(\mathbf{u}) \: dx
    \nonumber
\end{align}
The first integral can be rewritten as
\begin{align}
    &\frac{\partial}{\partial t} \int_{-\Delta_x/2}^{+\Delta_x/2}\mathbf{u}(x,y,t) \: dx
    = -\dot{s}[\mathbf{u}(s+,y,t)-\mathbf{u}(s-,y,t)] + \int_{-\Delta_x/2}^{+\Delta_x/2}\frac{\partial}{\partial t} \mathbf{u}(x,y,t^n) \:
    dx \nonumber
\end{align}
where $s+$ and $s-$ denote one-sided limits at the discontinuity
line. Thus if $\eqref{eq:Euler-int100}$ vanishes for $\Delta_x
\rightarrow 0$ we obtain the jump conditions $\eqref{eq:rankine}$
from the last two equations. The integral does not vanish, if
density perturbations can grow infinitely fast in time.
\newline
\newline
Viscosity and heat conduction in the Navier-Stokes equations
strongly affects steep gradients and perturbations. In reality a
contact surface cannot be maintained for an appreciable length of
time; (viscosity and) heat conduction between the permanently
adjacent particles on either side of the discontinuity would soon
make the idealized assumption unrealistic. While gas particles
crossing a shock front are exposed to heat conduction for only a
very short time, those that remain adjacent on either side of a
contact surface are exposed to heat conduction all the time. Hence a
contact layer will gradually fade out; see [CF1948]. Therefore some
odd behavior of the solution must be expected for a mathematical
idealized contact discontinuity.
\\
\section{Viscosity and Heat Conduction}
It is well known that viscosity and heat conduction can attenuate
small scale oscillations. In this section we study the effect of
heat conduction and viscosity on the solutions
$\eqref{eq:shea-flow-label}$ derived in the previous section.
Including viscosity in our considerations we obtain from the
Navier-Stokes equation for a divergence free flow
\begin{subequations}
\begin{align}
 &\frac{\partial}{\partial t}\rho + u \frac{\partial}{\partial
 x}\rho  + v \frac{\partial}{\partial y}\rho  = 0 \label{eq:density-000a} \\
 & \frac{\partial}{\partial t} u + u \frac{\partial}{\partial x} u + v \frac{\partial}{\partial y} u
 =  -\frac{1}{\rho}\frac{\partial}{\partial x}p
 + \frac{\mu}{\rho}(\frac{\partial^2}{\partial x^2}u + \frac{\partial^2}{\partial y^2}u)
 \label{eq:density-000b} \\
 &\frac{\partial}{\partial t} v + u \frac{\partial}{\partial x} v + v \frac{\partial}{\partial y} v
 = -\frac{1}{\rho}\frac{\partial}{\partial y}p
 + \frac{\mu}{\rho} (\frac{\partial^2}{\partial x^2}v + \frac{\partial^2}{\partial y^2}v)
  \label{eq:density-000c} \\
 &c_v \rho [\frac{\partial}{\partial t}T + u \frac{\partial}{\partial
    x}T + v \frac{\partial}{\partial y}T]
    = k (\frac{\partial^2}{\partial x^2}T + \frac{\partial^2}{\partial
    y^2}T)  \\
    &+2\mu (\frac{\partial}{\partial x}u )^2  + 2\mu (\frac{\partial}{\partial y}v
    )^2 + \mu (\frac{\partial}{\partial x}v + \frac{\partial}{\partial y}u )^2  \nonumber \\
    &\frac{\partial}{\partial x}u + \frac{\partial}{\partial y}v = 0
   \label{eq:density-000e}
\end{align} \label{eq:main-10000}
\end{subequations}
\newline
see e.g. [ATP1984;Section 5]. $T$ is the temperature given by
$\eqref{eq:eqof}$. For a constant x-component of the velocity
$u=u_0$ we obtain from the divergence condition
$\eqref{eq:density-000e}$, that the tangential component $v$ does
not depend on $y$ and equations $ \eqref{eq:main-10000}$ reduce to
\begin{subequations}
\begin{align}
 &\frac{\partial}{\partial t}\rho + u_0\frac{\partial}{\partial
 x}\rho + v \frac{\partial}{\partial y} \rho = 0
 \label{eq:NV-3a} \\
 &\frac{\partial}{\partial t}v + u_0 \frac{\partial}{\partial x}v
 =  -\frac{1}{\rho}\frac{\partial}{\partial y}p + \frac{\mu}{\rho} \frac{\partial^2}{\partial x^2}v
 \label{eq:NV-3c} \\
 & \frac{\partial}{\partial t}p + v \frac{\partial}{\partial y}p  = (\gamma -1)\mu (\frac{\partial}{\partial x}v)^2
  +  (\gamma -1)k (\frac{\partial^2}{\partial x^2}T + \frac{\partial^2}{\partial y^2}T)
  \label{eq:NV-3ee} \\
      &\frac{\partial}{\partial x}u + \frac{\partial}{\partial y}v = 0
   \label{eq:NV-3e}
\end{align}  \label{eq:NV-3-label}
\end{subequations} \newline
where we used $\eqref{eq:eqof}$ and $\eqref{eq:gamme-R-relation}$
and assumed that the viscosity $\mu$ and the coefficient of thermal
conductivity $k$ are positive constants. We obtain for a density of
the form $\eqref{eq:prop1-rho}$ with $u_r=u_0$ and
$v^{\varepsilon}=v$ for $x < u_0 t$ with $\eqref{eq:NV-3c}$
\begin{align}
&\frac{\partial}{\partial t}\rho + u_0 \frac{\partial}{\partial
    x}\rho + v \frac{\partial}{\partial y}\rho \nonumber \\
    &=-[\frac{\partial}{\partial t}v + u_0 \frac{\partial}{\partial
    x}v + \frac{\partial}{\partial y}v]t \frac{\partial}{\partial y}\rho_0
    \nonumber \\
    &=[\frac{1}{\rho}\frac{\partial}{\partial y}p - \frac{\mu}{\rho} \frac{\partial^2}{\partial x^2}v]t \frac{\partial}{\partial y}\rho_0
    \nonumber
\end{align}
For $t > 0$ $\eqref{eq:NV-3a}$ requires
\begin{align}
    &\frac{\partial}{\partial y}p - \mu \frac{\partial^2}{\partial x^2}v=0 \;\;\;\;\mathrm{or}\;\;\;\;\frac{\partial \rho_0}{\partial \eta} = 0 \label{eq:0-diff}
\end{align}
Since the x-component of the velocity is constant, the pressure
cannot depend on $x$. Since $v$ only depends on the spatial
coordinate $x$, this condition can only be satisfied if
\begin{align}
    &\frac{\partial}{\partial y}p = \mu \frac{\partial^2}{\partial x^2}v = f(t)\;\;\;\;\mathrm{or}\;\;\;\;\frac{\partial \rho_0}{\partial \eta} = 0 \label{eq:0-diff-1}
\end{align}
for a function $f(t)$. Therefore either the density does not depend
on $y$ or $v$ is a quadratic function in $x$; i.e. for a spatial
oscillatory tangential velocity $v$, the density $\rho$ cannot vary
in the y-direction. If viscosity is taken into account, perturbed
solutions of the form $\eqref{eq:shea-flow-label}$ with an
oscillatory tangential velocity $v_0^{\varepsilon}$ are excluded. In
a constant pressure flow any density perturbation is equivalent to
an entropy perturbation. Since $\eqref{eq:shea-flow-label}$ is for
$x\neq u_r t$ a smooth solutions of the Euler equations, these
solutions cannot be excluded through the established entropy
conditions for weak solutions. Thus we may regard
$\eqref{eq:0-diff-1}$ as an additional "plane wave condition" for
the Euler equations. A density perturbation of the form
$\eqref{eq:shea-flow-label}$ which violates $\eqref{eq:0-diff-1}$
will be denoted as an "inviscid entropy perturbation".
\newline
For a constant pressure $p=p_0$ equations $\eqref{eq:main-10000}$
reduce to
\begin{subequations}
\begin{align}
 &\frac{\partial}{\partial t}\rho + u \frac{\partial}{\partial
 x}\rho  + v \frac{\partial}{\partial y}\rho  = 0 \label{eq:density-p000a} \\
 & \frac{\partial}{\partial t} u + u \frac{\partial}{\partial x} u + v \frac{\partial}{\partial y} u
 = \frac{\mu}{\rho}(\frac{\partial^2}{\partial x^2}u + \frac{\partial^2}{\partial y^2}u)
 \label{eq:density-p000b} \\
 &\frac{\partial}{\partial t} v + u \frac{\partial}{\partial x} v + v \frac{\partial}{\partial y} v
 = \frac{\mu}{\rho} (\frac{\partial^2}{\partial x^2}v + \frac{\partial^2}{\partial y^2}v)
  \label{eq:density-p000c} \\
 & \frac{\partial^2}{\partial x^2}T + \frac{\partial^2}{\partial
    y^2}T = -2 \frac{\mu}{k} (\frac{\partial}{\partial x}u )^2  - 2 \frac{\mu}{k} (\frac{\partial}{\partial y}v
    )^2 - \frac{\mu}{k} (\frac{\partial}{\partial x}v + \frac{\partial}{\partial y}u )^2  \label{eq:density-p000d} \\
    &\frac{\partial}{\partial x}u + \frac{\partial}{\partial y}v = 0\label{eq:density-p000e}
\end{align} \label{eq:pmain-10000}
\end{subequations}
\newline
where we used $\eqref{eq:eqof}$ to rewrite the temperature as a
function of the pressure and the density and used
$\eqref{eq:density-p000a}$ to obtain $\eqref{eq:density-p000d}$. For
a given velocity field equation $\eqref{eq:density-p000d}$ is a
Poisson equation for the temperature; where the ratio of the
viscosity and thermal conductivity in $\eqref{eq:density-p000d}$ can
be computed from the Prandtl number $Pr$ and the specific heat at
constant pressure $c_p$ by
\begin{align}
    \frac{\mu}{k}=\frac{Pr}{c_p}
\end{align}
The ratio $Pr / c_p$ is approximately constant for most gases. For
air at standard conditions $Pr/c_p\simeq0,00072$.
\newline The temperature and constant
pressure defines the density via the equation of state; i.e.
\begin{align}
    \rho =\frac{p}{R} \frac{1}{T} \label{eq:T-rho-relation}
\end{align}
For a non vanishing coefficient of thermal conductivity $k$
$\eqref{eq:density-p000d}$ imposes a regularity condition on the
density, which is not present in an flow governed by the Euler
equations $\eqref{eq:Euler-diff}$; e.g. using the Weyl Lemma
[WA1994;Section 9] we obtain from $\eqref{eq:density-p000d}$ for a
velocity field which is constant outside a  bounded set, that the
density is as smooth as the velocity field, whereas the Euler
equations admits contact discontinuities for such a velocity field;
see Proposition 1.
\newline A similiar regularity condition holds for the pressure in
the incompressible Navier-Stokes equations. For a constant density
$\rho=\rho_0$ the equations $\eqref{eq:main-10000}$ reduce to
\begin{subequations}
\begin{align}
 & \frac{\partial}{\partial t} u + u \frac{\partial}{\partial x} u + v \frac{\partial}{\partial y} u
 =  -\frac{1}{\rho_0}\frac{\partial}{\partial x}p
 + \frac{\mu}{\rho_0}(\frac{\partial^2}{\partial x^2}u + \frac{\partial^2}{\partial y^2}u)
 \label{eq:density-rho000b} \\
 &\frac{\partial}{\partial t} v + u \frac{\partial}{\partial x} v + v \frac{\partial}{\partial y} v
 = -\frac{1}{\rho_0}\frac{\partial}{\partial y}p
 + \frac{\mu}{\rho_0} (\frac{\partial^2}{\partial x^2}v + \frac{\partial^2}{\partial y^2}v)
  \label{eq:density-rho000c} \\
 &\frac{\partial}{\partial t}p + u \frac{\partial}{\partial
    x}p + v \frac{\partial}{\partial y}p
    = \frac{k}{\rho_0 c_v} (\frac{\partial^2}{\partial x^2}p + \frac{\partial^2}{\partial
    y^2}p)  \label{eq:density-rho000d} \\
    &+ 2 (\gamma -1) \mu (\frac{\partial}{\partial x}u )^2  + 2 (\gamma -1)\mu (\frac{\partial}{\partial y}v
    )^2 + (\gamma -1) \mu (\frac{\partial}{\partial x}v + \frac{\partial}{\partial y}u )^2  \nonumber \\
    &\frac{\partial}{\partial x}u + \frac{\partial}{\partial y}v = 0
   \label{eq:density-rho000e}
\end{align} \label{eq:rhomain-10000}
\end{subequations}
where we used $\eqref{eq:T-rho-relation}$ and
$\eqref{eq:gamme-R-relation}$ to express the temperature as a
function of the density and pressure. Equation
$\eqref{eq:density-rho000b}$, $\eqref{eq:density-rho000c}$ and
$\eqref{eq:density-rho000e}$ constitute the incompressible
Navier-Stokes equations. The pressure in these equations can be
computed from the velocity field through a Poisson equation
\begin{align}
    \frac{\partial^2}{\partial x^2}p + \frac{\partial^2}{\partial
    y^2}p = - 2 \rho_0 (\frac{\partial}{\partial x}v \frac{\partial}{\partial y}u
    - \frac{\partial}{\partial x}u \frac{\partial}{\partial y}v)
   \label{eq:density-rho0001e}
\end{align}
see e.g. [KL1989; Section 9.1.3]. Again using the Weyl Lemma we
obtain from \eqref{eq:density-rho0001e} for a velocity field which
is constant outside a  bounded set, that the pressure is as smooth
as the velocity field. Inserting $\eqref{eq:density-rho0001e}$ into
$\eqref{eq:density-rho000d}$ we obtain for a given smooth velocity
field a inhomogeneous advection equation for the pressure. This
equation defines the time evolution of the pressure and is generally
neglected in an incompressible flow. Consider a constant pressure,
constant density flow with a constant x-component of the velocity
$u=u_r$. We obtain from $\eqref{eq:density-rho0001e}$ for a constant
pressure, that the functional determinate of $(u,v)$ vanishes.
Therefore $v=\Phi(u)$ for some smooth function $\Phi$ and we obtain
that the tangential velocity is constant as well. From proposition 1
follows that $(u_r,v)$ with $v=\Phi(x-u_r t)$ is a velocity field
for the Euler equations with a constant pressure and constant
density. These solutions are excluded for a finite viscosity $\mu$
and coefficient of thermal conductivity $k$ in the divergence free
Navier-Stokes equations $\eqref{eq:main-10000}$.
\newline
\newline
Numerical methods for the Euler equations employ an artificial
numerical viscosity model to resolve discontinuities. If this
numerical viscosity model is not properly related to the physical
viscosity and thermal conductivity, numerical artefact´s may be
introduced into a approximate solution, which are not related to a
real viscous, heat conducting flow.
\\
\section{Instability of Godunov´s method}
In a first order finite volume methods the assumption is generally
made, that the solution is constant inside a cell at a time-level
$t=t^n$; e.g.
\begin{align}
    \mathbf{\mathbf{\bar{u}}}(\vec{x},t):= \mathbf{\mathbf{u}}_{i,j} \;\;\mathrm{for}
    \;\;\vec{x}\;\epsilon \;[x_{i-1/2},x_{i+1/2}]\times[y_{i-1/2},y_{i+1/2}]
\end{align}
and it is assumed that discontinuities are moved to the cell
boundary. We assume for simplicity that
$D_{i,j}=[x_{i-1/2},x_{i+1/2}]\times[y_{i-1/2},y_{i+1/2}]$ is a
rectangle defined through a constant cartesian grid
${x_{i}=i\Delta_x}$ and ${y_{j}=j\Delta_y}$ and denote by
\begin{align}
    \mathbf{\bar{u}}_{i,j}(t)= \frac{1}{V(D_{i,j})}\int_{D_{i,j}}^{}\mathbf{u}(\xi,\eta,t) \: d\eta \; d\xi  \label{eq:cell-ave1}
\end{align}
the cell-average, where $V(D_{i,j})$ is the Volume of $D_{i,j}$.
Finite-Volume Godunov-type methods are derived from the integral
form of the conservation law. The integral form
$\eqref{eq:Euler-int}$ can be rewritten as
\begin{align}
    & \frac{d}{dt}\mathbf{\bar{u}}_{i,j}(t)= -\frac{1}{V(D_{i,j})} \int_{\partial D_{i,j}}^{ } \vec{n}(\xi) \cdot\: \vec{F}(\mathbf{u}^R(\vec{x}(\xi),t+;\vec{n}(\xi))) \: d\xi \label{eq:Euler-int2}
\end{align}
This equations says that we can evolve the cell averages in time, by
solving one-dimensional Riemann problems at a the cell boundary at
time $t=t^n$ and then solve the system of ordinary differential
equation $\eqref{eq:Euler-int2}$ to obtain the cell average at time
$t=t^n+\tau$, $\tau > 0$. Taking for granted that the solution of
the Riemann problem at the cell interface can be locally advanced in
time; i.e. we require that
\begin{align}
    \frac{\partial}{\partial t}\mathbf{u}^R(\vec{x}(\xi),t;\vec{n}(\xi))|_{t^=t^n+} \label{eq:dt-cel-bou}
\end{align}
should exist at the cell-boundary $\partial D_{i,j}$.
\newline
In Godunov´s method a plane-wave Riemann problem is solved at the
cell boundary, for example at the cell-boundary $(x_{i+1/2},y_j)$ at
time $t=t^n$ in the x-direction. Denote by
$\mathbf{u}^R(x_{i+1/2},y_j,t)$ the solution of the plane-wave
Riemann problem in the x-direction at the cell boundary, then
$\mathbf{u}^n_{i+1/2,j}:=\lim_{t\rightarrow t^n+}
\mathbf{u}^R(x_{i+1/2},y_j,t)$ is computed and used to evaluate the
physical flux function.
\newline
In the last two sections we saw that a constant pressure region
which a contact discontinuity is unstable under perturbations. In
this section we analyze Godunov´s method for this critical region.
Assume that the velocity $u$ and the pressure $p$ is constant, then
the solution of the Riemann problem in the x-direction consists of
contact discontinuity and we obtain
\begin{align}
    \mathbf{u}^R(x_{i+1/2},y_j,t^n)= \label{eq:x-rieman-contact}
    \begin{cases}
        \mathbf{u}_{i+1,j} & \text{for  } u < 0 \\
        \mathbf{u}_{i,j}   & \text{for  } 0 < u
    \end{cases}
\end{align}
If $v$ depends only on the spatial coordinate x and $p$ is constant,
then the solution in the y-direction is also a contact discontinuity
and we obtain
\begin{align}
    \mathbf{u}^R(x_i,y_{j+1/2},t^n)= \label{eq:y-rieman-contact}
    \begin{cases}
        \mathbf{\bar{u}}_{i,j+1} & \text{for  } v_i < 0 \\
        \mathbf{\bar{u}}_{i,j}   & \text{for  } 0 < v_i
    \end{cases}
\end{align}
Therefore for $u < 0$ and $v_i < 0$ we obtain
\begin{align}
    \mathbf{\bar{u}}^{n+1}_{i,j}&=\mathbf{\bar{u}}^{n}_{i,j}
    - \frac{\tau}{\Delta_x} [\mathbf{f}(\mathbf{\bar{u}}^n_{i+1,j}) - \mathbf{f}(\mathbf{\bar{u}}^n_{i,j})]
    - \frac{\tau}{\Delta_y} [\mathbf{g}(\mathbf{\bar{u}}^n_{i,j+1}) - \mathbf{g}(\mathbf{\bar{u}}^n_{i,j})]
    \label{eq:godunov-1}
\end{align}
which is equivalent to the four difference equations
\begin{align}
    \bar{\rho}^{n+1}_{i,j} &= \bar{\rho}^{n}_{i,j} - \frac{\tau}{\Delta_x} u (\bar{\rho}^{n}_{i+1,j} - \bar{\rho}^{n}_{i,j})
                        - \frac{\tau}{\Delta_y} \bar{v}_{i} (\bar{\rho}^{n}_{i,j+1} - \bar{\rho}^{n}_{i,j}) \nonumber \\
    u \bar{\rho}^{n+1}_{i,j} &= u[\bar{\rho}^{n}_{i,j} - \frac{\tau}{\Delta_x} u (\bar{\rho}^{n}_{i+1,j} - \bar{\rho}^{n}_{i,j})
                        - \frac{\tau}{\Delta_y} \bar{v}_{i} (\bar{\rho}^{n}_{i,j+1} - \bar{\rho}^{n}_{i,j})] \nonumber \\
    (\bar{\rho} \bar{v})^{n+1}_{i,j} &= (\bar{\rho} \bar{v})^{n}_{i,j} - \frac{\tau}{\Delta_x} u ((\bar{\rho} \bar{v})^{n}_{i+1,j} - (\bar{\rho} \bar{v})^{n}_{i,j})
                        - \frac{\tau}{\Delta_y} \bar{v}_{i} ((\bar{\rho} \bar{v})^{n}_{i,j+1} - (\bar{\rho} \bar{v})^{n}_{i,j}) \nonumber \\
    \bar{E}^{n+1}_{i,j} &= \bar{E}^{n}_{i,j} - \frac{\tau}{\Delta_x} u (\bar{E}^{n}_{i+1,j} - \bar{E}^{n}_{i,j})
                        - \frac{\tau}{\Delta_y} \bar{v}_{i} (\bar{E}^{n}_{i,j+1} -
                        \bar{E}^{n}_{i,j})   \nonumber
\end{align}
Since the second equation is just the first equation multiplied with
the constant velocity $u$, this can be reduced to three difference
equations
\begin{align}
    \bar{\rho}^{n+1}_{i,j} &= \bar{\rho}^{n}_{i,j} - \frac{\tau}{\Delta_x} u (\bar{\rho}^{n}_{i+1,j} - \bar{\rho}^{n}_{i,j})
                        - \frac{\tau}{\Delta_y} \bar{v}_{i} (\bar{\rho}^{n}_{i,j+1} - \bar{\rho}^{n}_{i,j}) \nonumber \\
   (\bar{\rho} \bar{v})^{n+1}_{i,j} &= (\bar{\rho} \bar{v})^{n}_{i,j} - \frac{\tau}{\Delta_x} u ((\bar{\rho} \bar{v})^{n}_{i+1,j} - (\bar{\rho} \bar{v})^{n}_{i,j})
                        - \frac{\tau}{\Delta_y} \bar{v}_{i} ((\bar{\rho} \bar{v})^{n}_{i,j+1} - (\bar{\rho} \bar{v})^{n}_{i,j}) \nonumber \\
    \bar{E}^{n+1}_{i,j} &= \bar{E}^{n}_{i,j} - \frac{\tau}{\Delta_x} u (\bar{E}^{n}_{i+1,j} - \bar{E}^{n}_{i,j})
                        - \frac{\tau}{\Delta_y} \bar{v}_{i} (\bar{E}^{n}_{i,j+1} -
                        \bar{E}^{n}_{i,j})   \nonumber
\end{align}
This is a discrete approximation to
\begin{subequations}
\begin{align}
    &\frac{\partial}{\partial t}\bar{\rho} + u \frac{\partial}{\partial x}\bar{\rho} + \bar{v} \frac{\partial}{\partial y}\bar{\rho} = 0 \label{eq:shear-flow-1000} \\
    &\frac{\partial}{\partial t}(\bar{\rho} \bar{v}) + u \frac{\partial}{\partial x}(\bar{\rho} \bar{v}) + \bar{v} \frac{\partial}{\partial y}(\bar{\rho} \bar{v}) = 0 \label{eq:shear-flow-1001} \\
    &\frac{\partial}{\partial t}\bar{E} + u \frac{\partial}{\partial x}\bar{E} + \bar{v} \frac{\partial}{\partial y}\bar{E} = 0 \label{eq:shear-flow-1002}
\end{align} \label{eq:shear-flow-2000}
\end{subequations}
Rewriting the second equation as
\begin{align}
        &\bar{\rho}[\frac{\partial}{\partial t}\bar{v} + u \frac{\partial}{\partial x}\bar{v}]
        +\bar{v}[\frac{\partial}{\partial t}\bar{\rho} + u \frac{\partial}{\partial x}\bar{\rho} + \bar{v} \frac{\partial}{\partial y}\bar{\rho}] = 0 \nonumber
\end{align}
we obtain that $\eqref{eq:shear-flow-2000}$, is equivalent to
$\eqref{eq:shear-flow}$ in this flow region. This holds also for $u
< 0$, $v \geq 0$, $u \geq 0$, $v < 0$ and $u \geq 0$, $v \geq 0$.
Therefore Godunov´s scheme is a discrete approximation to
$\eqref{eq:shear-flow}$ in a regions with a constant pressure $p$,
constant x-component of the velocity $u$ and a y-component of the
velocity $v$, which does only depend on the spatial coordinate $x$.
\newline
For a stationary shock wave at $x_{1/2}=0$ with a shock line
perpendicular to the x-axis, the x-component of the flux function of
Godunov´s scheme satisfies
\begin{align}
     \mathbf{f}^{G}_{1/2,j}=\mathbf{f}(\mathbf{\bar{u}}^n_{1,j}) =
     \mathbf{f}(\mathbf{\bar{u}}^n_{0,j})
\end{align}
i.e. one-dimensional stationary shocks are resolved exactly.
Therefore Godunov´s method is a discrete approximation to
$\eqref{eq:shear-flow}$ behind the shock with boundary conditions
\begin{align}
    \bar{\rho}(0,t)=const. \nonumber \\
    \bar{v}(0,t)=const. \nonumber
\end{align}
where $u_r$ in $\eqref{eq:shear-flow}$ is the constant normal
component of the velocity behind the shock. If the shock is near
stationary these boundary conditions introduce perturbation. Thus on
a sufficiently fine grid entropy perturbations of the form
$\eqref{eq:shea-flow-label}$ and acoustic perturbations of the form
$\eqref{eq:isenf-per-rho}$ are resolved by Godunov´s method behind a
near stationary shock. In the last section we saw that inviscid
entropy perturbations do not relate to a real viscous flow. Flow
structures not related to a real viscous, heat conducting flow may
therefore appear in numerical solutions of Godunov´s method.
\newline
\\
Numerical examples for a plane shock wave aligned with the grid,
which is moving down a duct are given in [QUI 1994; Figure 5]. At
the grid center line, a small perturbation is introduced in the
computation. Downstream of the shock an unstable density profile
develops, which over time leads to an unstable numerical shock
front. If we associate the center of the duct with the x-axis, then
the perturbation depend behind the shock front on $y$. This
numerical example reflects the situation discussed analytically at
the beginning of section 2.
\newline
Roe´s method [Roe1980] is more prone to generate these perturbations
then Godunov´s method, due to the fact that a rarefaction wave is
replaced by a rarefaction shock. Therefore noise from the
rarefaction shock, can also lead to an unstable growth of the
density, in addition to the noise from the shock.
\newline
An example from an aerodynamic simulation, which results in
incorrect numerical results is given in [PEIM1988], for a bow shock
over a blunt body placed in a high Mach number flow. Along the
stagnation line the bow shock is approximately aligned with grid
used for the calculation. A perturbation of the shock profile is
given through the curvature of the shock. At the stagnation point we
have approximately a plane-wave near stationary shock wave, with a
disturbed shock profile. This again is the situation discussed at
the beginning of section 2.
\\
Based on the numerical observation, that shock capturing methods
which try to capture contact discontinuities exactly, generally
suffer from failings, a link between the carbuncle phenomenon and
the resolution of the contact discontinuities was suggested by
Gressier and Moschetta [GRE1998].
\\
The dissipation model for Godunov´s scheme was studied by Xu
[XU1999] in a series of numerical experiments.  He concluded that
Godunov´s method gives accurate results in both unsteady shock
structure and boundary layer calculations, but that the absence of
dissipation in the gas evolution model in Godunov´s scheme amplifies
post-shock oscillations. He found that the numerical dissipation
model for Godunov´s method is in the multidimensional case
mesh-oriented and not consistent with the Navier-Stokes equations.
Xu also mentioned that it is well known that the inviscid Euler
equations cannot give a correct representation of the fluid motion
in the discontinuity flow region. Which is consistent with our
analysis regarding a contact discontinuity line (vortex sheet).
\newline Real weak shock fronts are transition layers of finite
width and the representation of weak shock fronts through a
discontinuity line in the Euler equations is a mathematical
approximation. A justification for this approximation was given by
Majda [MA1983; Proposition 3]. Based on a linear stability analysis
Majda found that planar compressive shock fronts in an ideal gas are
uniformly stable. In [LL1959; §87] it is furthermore noted that the
front depth of a non weak shocks is so small, that a transition
layer becomes meaningless. Also by no means obvious, it may be
safely assumed, that the mathematical representation of a single
smooth shock front in an ideal gas through a discontinuity line
resp. surface is reasonable. Since Majda´s stability analysis does
not apply to a contact discontinuity line/surface and we assumed in
this paper that shock fronts only generate perturbations without
assuming an unstable growth at the smooth shock front itself, this
result is consistent with the analysis in this paper.
\newline The source of small perturbations in Godunovs methods is
the displacement of the shock curves at the cell-boundary and the
nonlinear interaction of the dependent variables in the numerical
shock layer. A non stationary or stationary displaced shock wave is
approximated through a smeared profile with at least one
intermediate cell. Whenever the smeared shock profile changes,
perturbations are generated from the characteristic fields. If the
shock curve is not exactly a plane-wave in the x-direction, the
shock-curvature will introduce an $y$ dependence in the
perturbations. In regions of low (numerical) viscosity, the
instability under perturbation of the constant pressure region
behind the shock  will result in an unstable flow - if a contact
discontinuity is present.
\newline For the case reported by Perry
and Imlay [PEIM1988] entropy perturbations are generated a the
shock. Since according to the proposition in section 2, entropy
perturbations can only exits if the tangential velocity does not
depend on $y$, these perturbations are initially restricted to the
symmetry axis. If the flow is near stationary, cell-averaging of the
form $\eqref{eq:cell-ave1}$ in the projection state in Godunov´s
method leads to density perturbation in these cells. At the boundary
of these cells Godunov´s scheme then resolves a contact
discontinuity, which results in an unstable flow according to the
analysis in section 2 and 4. The unstable flow behind the shock
interacts with the near stationary shock front resulting in
carbuncle structures. Since the interaction destroy´s also the
smoothness of the shock front itself, the stability analysis of
Majda no longer applies.
\newline
For the case reported by Quirk [QUI 1994; Figure 5] perturbation are
introduced externally at the grid center line. The grid center line
corresponds to the symmetry line in results of [PEIM1988]. In [QUI
1994; Figure 5a] perturbations are first visible behind the plane
shock, but no carbuncle is visible at this stage. Only after the
shock has propagated further down the duct the carbuncle becomes
visible in front of the shock.
\newline
If discontinuous density or velocity perturbation are artificially
introduced in a background flow then according to the analysis in
section 2, unstable structures can even be generated in numerical
solutions of the isentropic Euler equations. This was demonstrated
by Elling [VE2006] for Godunov´s method. In his case a steady plane
shock line parallel to the y-axis is perturbed through a
one-cell-high filament along the x-axis in front of the shock. In
the filament the normal component of the velocity $u=\bar{u}$ is set
to zero. The flow in the filament is therefore a constant pressure
flow governed by $\eqref{eq:pconst}$ with a constant normal velocity
component $v=\bar{v}$. From proposition 1 follows (with the
tangential and normal velocity interchanged) that
$\bar{u}=\bar{u}(y,t)$. The linear stability analysis results again
in an advection equation with a singular source term of the form
 \begin{align}
 & \frac{\partial}{\partial t}\tilde{u} + \bar{v} \frac{\partial}{\partial y}\tilde{u}
   = - \tilde{v} (\bar{u}_r - \bar{u}_l)\delta(y - \bar{v} t) \label{eq:epw-simple-carbuncle}
 \end{align}
This shows that the flow in the filament is unstable. Interactions
of this unstable flow with the shock wave results in a carbuncle
instability. Elling noted that the carbuncle instability can also be
observed for the local Lax-Freidrich/Rusanov and Osher-Solomon
schemes. His conjecture is that carbuncles can be related to a
special class of non-physical entropy solutions for the continuum
equations, which is supported by our analysis.
\\
A discussion of the carbuncle instability for several numerical
methods, can be found in [DMG;2004] and [Roe2007]. The focus in
these papers is on the perturbations generated through a numerical
shock profile. The cause for the carbuncle instability is often
related to the numerical shock profile itself and the instability
also denoted a as a shock instability. However, in this paper it is
shown, that the instability may be related to the central constant
pressure region in the plane wave Riemann problem. This is an
intrinsic instability of the Euler equations and not related to a
particular numerical method. The numerical shock profile manifests
this instability in regions of low numerical dissipation through the
generation of perturbations. This manifestation of the instability,
is problem and method dependent.
\newline
%
\\
\section{Summary and Conclusion}
In this paper it is proved that a constant pressure flow region,
governed by the hyperbolic conservations laws for an ideal gas, is
unstable under perturbations if a discontinuity is present. The
instability is immanent to the linear degenerate fields in the
multidimensional hyperbolic conservation laws, with a double
eigenvalue.
\newline
The mathematical idealized assumption of a zero thickness contact
discontinuity in the plane wave Riemann problem is unrealistic for
real flows. If infinitesimal viscosity and thermal conductivity are
taken into account, inviscid entropy perturbations are excluded for
an oscillatory tangential velocity in a perturbed plane wave Riemann
problem. Since smooth entropy disturbances in a constant pressure
flow are transported with the fluid, they cannot be excluded through
the established entropy conditions for weak solutions of hyperbolic
conservations laws. Additional conditions are required to guarantee
that a solution of the plane wave Riemann problem is as a limit
solution of the Navier-Stokes equations for a vanishing viscosity
and thermal conductivity.
\newline
Immanent perturbations (acoustic and entropy) generated at a
perturbed shock line and the numerical impossibility to
differentiate exactly between admissible and non admissible entropy
perturbations in a high Reynols number (low viscous) flow, are
challenges for discontinuity resolving methods. Godunov´s method
closely relates to the physics of the Euler equations. Numerical
artefact´s observed in numerical solutions for Godunov´s method are
a numerical manifestation of the immanent instability in the Euler
equations. The HLLE scheme [EIN1988] on the other hand can be
regarded as a vortex sheet averaged approximation to the
Navier-Stokes equations for high Reynolds number flows.
\newpage
\section{References}
$\\$[ATP1984] D. A. Anderson, J. C. Tannehill, R. H. Pletcher,
"Computational Fluid Mechanics and Heat Transfer", HEMISPHERE
Publishing Corporation, 1984.$\\$ [CF1948] R. Courant and K. 0.
Friedrichs, "Supersonic Flow and Shock Waves", Interscience, New
York, 1948 $\\$[DMG;2004]M. Dumbser, J. M. Moschetta and J.
Gressier, "A Matrix stability analysis of the carbuncle phenomenon",
J. Comput. Phys., 197, 647-670 (2004)$\\$ [EIN1988] B. Einfeldt, "On
Godunov-type methods for gas dynamics", SIAM J. Numer. Anal., 25,
294-318 (1988). $\\$ [GRE1998]J. Gressier and J.-M. Moschetta
"Robustness versus Accuracy in Shock-Wave Computations",
International Journal of Numerical Methods in Fluids, July 1999.$\\$
[KL1989] H. O. Kreiss and J. Lorenz, "Initial-Boundary Value
Problems and the Navier-Stokes Equations", Academic Press 1989. $\\$
[LAX2006] P. D. Lax, "Hyperbolic Partial Differential Equations",
Courant Lecture Notes 14, 2006. $\\$[LEV2002] R. J. LeVeque, "Finite
Volume Methods for Hyperbolic Problems", Cambridge University Press,
Cambridge, 2002.$\\$ [LL1959] L. D. Landau and E. M. Lifchitz,
"Fluid Mechanics", Addison-Wesley Publishing, 1957.$\\$[MA1983] A.
Majda, "The stability of multi-dimensional shock fronts", Memories
of the AMS Number 275, 1983.$\\$ [PEIM1988] K.M. Perry and S. T
Imlay, "Blunt Body Flow Simulations", AIAA Paper, 88-2924,1988. $\\$
[QUI1994] J.J. Quirk, "A Contribution to the Great Riemann Solver
Debate", International Journal For Numerical Methods in Fluids, Vol.
18, 555-574 (1994) $\\$ [ROE2007] P. L. Roe, "An Evaluation of Euler
Fluxes for Hypersonic Flow Computations", in AIAA Computational
Fluid Dynamics Conference, Miami, FL, Jun. 25-28, 2007 $\\$ [SM1983]
J. Smoller, "Shock Waves and Reaction-Diffusion Equations",
Springer, 1983.$\\$[VE2006] V. Elling, "Carbuncles as self-similar
entropy solutions", arXiv:math/0609666v1, 2006$\\$[WA1994] W.
Walter, "Einführung in die Theorie der Distributionen", BI 1994 $\\$
[SM1983] J. Smoller, "Shock Waves and Reaction-Diffusion Equations",
Springer, 1983.$\\$ [XU1999] Kun Xu, "Gas Evolution Dynamics in
Godunov-type Schemes and Analysis of Numerical Shock Instability";
ICASE Report No. 99-6, 1999.
\end{document}